%
%
%
%
%
%
%
%
%
%
%
%
\scrollmode
\magnification=\magstep1
\parskip=\smallskipamount
\hoffset=1cm \hsize=12cm

\def\demo#1:{\par\medskip\noindent\it{#1}. \rm}
\def\ni{\noindent}               
\def\ll{\leftline}
\def\cl{\centerline}

\def\begin{\ll{}\vskip 5mm\nopagenumbers}  
\def\pn{\footline={\hss\tenrm\folio\hss}}   

%
%
\outer\def\beginsection#1\par{\bigskip
  \message{#1}\leftline{\bf\&#1}
  \nobreak\smallskip\vskip-\parskip\noindent}

%
%
\outer\def\proclaim#1:#2\par{\medbreak\vskip-\parskip
    \noindent{\bf#1.\enspace}{\sl#2}
  \ifdim\lastskip<\medskipamount \removelastskip\penalty55\medskip\fi}

\def\endpr{\hfill $\spadesuit$ \medskip}

%
%
%
%


%
%

\def\R{{\rm I\kern-0.2em R\kern0.2em \kern-0.2em}}
\def\N{{\rm I\kern-0.2em N\kern0.2em \kern-0.2em}}
\def\P{{\rm I\kern-0.2em P\kern0.2em \kern-0.2em}}
\def\B{{\rm I\kern-0.2em B\kern0.2em \kern-0.2em}}
\def\C{{\rm C\kern-.4em {\vrule height1.4ex width.08em depth-.04ex}\;}}
\def\CP{\C\P}
\def\RP{\R\P}

%
%
%
%

\def\cC{{\cal C}}

\def\cO{{\cal O}}

%
%
%
\def\a{\alpha}
\def\b{\beta}

\def\d{\delta}
\def\e{\epsilon}
\def\z{\zeta}

\def\l{\lambda}

\def\c{\chi}

\def\L{\Lambda}

%
%
%
%
\def\bar{\overline}              
\def\bs{\backslash}              
\def\di{\partial}                
\def\hra{\hookrightarrow}

%
%
\def\disc{\triangle}             

\def\c*{{\C}^*}

\def\wt{\widetilde}

%
%
\def\holo{holomorphic}                   
\def\nbd{neighborhood}                   
\def\psc{pseudoconvex}                   
\def\spsc{strongly\ pseudoconvex}        
\def\psh{plurisubharmonic}               
\def\spsh{strongly\ plurisubharmonic}
\def\tr{totally real}                    
\def\ss{\subset\!\subset}                

\def\iff{if and only if}


\begin
\cl{\bf  STEIN DOMAINS IN COMPLEX SURFACES}
\bigskip
\cl{Franc Forstneri\v c}

\vskip 10mm
\cl{\bf Contents}
\smallskip
\ll{Introduction \dotfill $\quad 1$}
\smallskip
\ll{1. Regular Stein neighborhoods of embedded surfaces \dotfill $\quad 3$}
\smallskip
\ll{2. Construction of regular Stein neighborhoods  \dotfill $\quad 7$}
\smallskip
\ll{3. Connected sums of embedded and immersed surfaces \dotfill $\,\,\, 12$}
\smallskip
\ll{4. Regular Stein neighborhoods of immersed surfaces \dotfill $\,\,\, 13$}
\smallskip
\ll{Appendix: Generalized adjunction inequalities \dotfill $\,\,\, 17$}
\smallskip
\ll{References  \dotfill $\,\,\, 20$}

\vskip 6mm

%
%
%
%
\cl{\bf Abstract}
Let $S$ be a  closed connected real surface and $\pi\colon S\to X$ a smooth embedding 
or immersion of $S$ into a complex surface $X$. We denote by $I(\pi)$ 
(resp.\ by $I_\pm(\pi)$ if $S$ is oriented) the number of complex points of 
$\pi(S)\subset X$ counted with algebraic multiplicities. Assuming that $I(\pi)\le 0$ 
(resp.\ $I_\pm(\pi)\le 0$ if $S$ is oriented) we prove that $\pi$ can be $\cC^0$ approximated 
by an isotopic immersion $\pi_1\colon S\to X$ whose image has a basis of open 
Stein \nbd s in $X$ which are homotopy equivalent to $\pi_1(S)$. We obtain precise 
results for surfaces in $\CP^2$ and find an immersed symplectic sphere in $\CP^2$ 
with a Stein \nbd.

\vskip 5mm
%
%
%
%
\cl{\bf Introduction}

A complex manifold $\Omega$ is called {\it Stein} if it is biholomorphic to a 
closed complex submanifold of a Euclidean space $\C^N$. Stein manifolds are of special
interest in complex analysis due to their rich function theory. Since every $n$-dimensional 
Stein manifold is homotopy equivalent to a CW-complex of real dimension 
at most $n$ [AF], it is natural to ask which such complexes can be realized 
by homologically nontrivial Stein domains in a given complex manifold. 

In this paper we consider this problem for CW complexes which are represented by 
smooth embeddings (or immersions) of closed connected real surfaces $S$ in complex 
surfaces $X$. Let $\pi\colon S\to X$ be such an immersion. 
An open \nbd\ $\Omega \subset X$ of $\pi(S)$ is called {\it regular} 
if it has a strong deformation retraction onto $\pi(S)$ (and hence $\Omega$ 
is homotopy equivalent to $\pi(S)$). In [F2] we proved that 
every closed real surface $S$ (not necessarily orientable) different 
from the two-sphere admits a smooth embedding in $\C^2$ with a basis of 
regular Stein \nbd s. Here we consider more general complex surfaces with 
emphasis on the complex projective plane $X=\CP^2$. 

Denote by $I(\pi) \in Z$ the number of complex points of $\pi$ counted 
with algebraic multiplicities (Sect.\ 1). If $S$ is oriented we denote by 
$I_\pm(\pi)$ the number if positive (resp.\ negative) complex points.
These numbers are sometimes called  {\it Lai indices} after [La],
although they appeared already in [CS] and [B]. Our main result is that, 
if $I_\pm(\pi) \le 0$ (resp.\ $I(\pi)\le 0$ when $S$ is unorientable) 
then $\pi$ can be $\cC^0$-approximated by a regularly homotopic immersion 
$\pi_1\colon S\to X$ whose image $S_1=\pi_1(S)\subset X$ has a basis of Stein 
\nbd s in $X$ homotopically equivalent to $S_1$ (Theorem 1.1). 
These \nbd s are sublevel sets of a smooth nonnegative weakly \psh\ function 
in a \nbd\ of $S_1$ which vanishes quadratically on $S_1$ and has no other 
critical points nearby (Theorem 2.2). We obtain precise results 
on the existence of such embeddings in the projective plane $\CP^2$ 
(Proposition 1.4 and Theorems 1.6 and 4.3). 

\pn

For oriented embedded real surfaces $S\subset X$ a converse to Theorem 1.1 is 
provided by gauge theory. It seems that the connection was first observed by
Lisca and Mati\'c [LM, Theorem 5.2] and S.\ Nemirovski [N, Theorem 9].
The condition $I_\pm(S) \le 0$ is equivalent to 
the {\it generalized adjunction inequality}
$$
	g(S)\ge 1+{1\over 2}(S^2+|c_1(X)\cdotp S|),                  \eqno(*)
$$
where $g(S)$ is the genus of $S$, $S^2=[S]^2\in Z$ is the self-intersection 
number of $[S]\in H_2(X;Z)$ in $X$ and $c_1(X) = c_1(\Lambda^2 TX)\in H^2(X;Z)$ 
is the first Chern class of $X$ (Sect.\ 1). The inequality (*) has been proved
in many interesting cases using gauge theory methods. The recent Seiberg-Witten 
approach caused a major revolution in this field and the literature 
is still growing rapidly. In the Appendix we summarize the results on (*) which are 
most relevant to us. In particular, (*) {\it holds if $X$ is a Stein surface and 
$S\subset X$ is any oriented, connected, embedded real surface except a homologically 
trivial two-sphere} (Theorem II in the Appendix). Hence the condition $I_\pm(S)\le 0$ 
is also necessary for the existence of a Stein \nbd\ of an oriented, 
embedded, homologically nontrivial surface. Nemirovski used this to solve 
the Vitushkin conjecture [N, Theorem 15] to the effect that there are no nonconstant
holomorphic functions in any \nbd\ of an embedded homologically nontrivial sphere 
in $\CP^2$. The corresponding generalized adjunction inequality for immersed 
oriented surfaces with simple double points is 
$$   g(S)+ \d_+ 
	\ge 1+ {1\over 2}
        \left( \pi(S)^2+|c_1(X)\cdotp \pi(S)| \right),     \eqno(**) 
$$
where $\d_+$ denotes the number of positive double points of $\pi$.
This holds if $\pi(S)$ has a Stein \nbd\ $\Omega\subset X$ and 
$[\pi(S)]\ne 0\in H_2(\Omega;Z)$ (Theorem III in the Appendix).
For $X=\CP^2$ we find embeddings and immersions which realize the lower
bounds in (*) resp.\ (**) (Theorem 4.3). 

Our examples, together with a symplectic approximation theorem of Gromov, 
give an immersed symplectic sphere $\pi\colon S\to \CP^2$ 
with a Stein \nbd\ (Corollary 4.4). A theorem of Ivashkovich and Shevchishin [IS] 
shows that such $\pi$ must have some negative double points since otherwise the 
envelope of meromorphy of any \nbd\ of $\pi(S)$ would contain a 
rational curve which would contradict the existence of a Stein \nbd. 
It is not clear what is the minimal number of simple double points of 
immersed symplectic spheres in $\CP^2$ which admit a Stein \nbd.

One can find Stein domains of more general homotopy type in a given complex 
surface by attaching handles as in [E1, Go] and using the Kirby calculus [GS]. 
We shall not pursue these matters here. 

\beginsection 1. Regular Stein \nbd s of embedded surfaces.

Let $S$ be a closed ($=$ compact without boundary) connected real surface 
smoothly embedded in a complex surface $X$ ($=$ a complex two dimensional
manifold). A point $p\in S$ is said to be {\it complex} 
if the tangent plane $T_p S$ is a complex line in 
$T_p X$; the other points are called {\it totally real}.
In suitable local \holo\ coordinates $(z,w)$ on $X$ near a complex point 
$p\in S$ the surface $S$ is the graph $w=f(z)$ of a smooth
function over a domain in $\C$, and $(z,f(z))\in S$ 
is a complex point of $S$ \iff\  
${\di f\over \di \bar z}(z) = {1\over 2}({\di f\over \di x}+i{\di f\over \di y})(z)= 0$. 
If $p=0$ is an isolated complex point of $S$ then for sufficiently small $\d>0$ we have
${\di f \over \di\bar z}(z)\ne 0$ for $0<|z|\le \d$, and we define the {\it index} 
$I(p;S)$ as the winding number of ${\di f \over \di\bar z}(z)$ on the circle 
$|z|=\d$. The index is independent of the choice of local coordinates 
and is related to the {\it Maslov index} [F1, F2]. A \tr\ point has index zero. 

For a generic choice of the embedding $S\subset X$ there are at most finitely many 
complex points in $S$ and one defines the index of the embedding by
$ 
	I(S)= \sum_{p\in S} I(p;S).  
$
If $S$ is oriented we divide its complex points into {\it positive}
resp.\ {\it negative}, depending on whether the orientation on 
$T_p S$ induced by the complex structure on $X$ agrees or disagrees
with the given orientation on $S$. In this case one defines 
$I_\pm(S) \in Z$ as the sum of indices over all positive 
(resp.\ negative) complex points of $S$. 
These notions clearly extend to immersions $\pi\colon S\to X$ 
in which case we write $I(\pi)=\sum_{p\in S} I(p;\pi)$ and similarly
for $I_\pm(\pi)$.

\proclaim Definition: 
A family of open sets $\{\Omega_\e \colon \e\in (0,\e_0) \}$ in 
$X$ is said to be a {\it regular basis of \nbd s of $S \subset X$} 
if for each $\e\in(0,\e_0)$ we have $\Omega_\e=\cup_{t<\e} \Omega_t$, 
$\bar \Omega_\e=\cap_{t>\e}\Omega_t$, and $S=\bigcap_{0<t<\e_0} \Omega_t$ 
is a strong deformation retraction of $\Omega_\e$. In particular,
every $\Omega_\e$ is homotopy equivalent to $S$.

The following is our main result. 

\proclaim 1.1 Theorem: Let $S$ be a closed, connected real surface 
smoothly embedded in a complex surface $X$. If $S$ is oriented 
and $I_\pm(S)\le 0$, or if $S$ is unorientable and $I(S)\le 0$, then
$S$ can be $\cC^0$-approximated by an isotopic embedding $S'\subset X$ with 
a regular basis of smoothly bounded Stein \nbd s. 
The analogous result holds for immersions.

Theorem 1.1 is proved in section 2. The special case $X=\C^2$
was considered in [F2] where it was proved that {\it every closed real surface
except the two-sphere admits an embedding in $\C^2$ with a regular
Stein \nbd\ basis}.

%
%
Before proceeding we recall the {\it index formulas} expressing the indices 
$I(\pi)$ (resp.\ $I_\pm(\pi)$) by invariants of the immersion. 
Isolated complex points of real surfaces in complex surfaces were first 
investigated by Chern and Spanier [CS]. Bishop [B] classified them  
into {\it elliptic\/} (these have index $+1$), 
{\it hyperbolic} (with index $-1$), and {\it parabolic} 
(these are degenerate and may have any index). If $S\subset X$
has $e$ elliptic points and $h$ hyperbolic points then $I(S)=e-h$
(this happens for a generic $S$). Bishop [B] proved that for every oriented 
embedded surface $S\subset \C^2$ we have $I(S)= \chi(S)=2-2g(S)$ 
and $I_\pm(S)={1\over 2}I(S)=1-g(S)$. 
In general we have the following ([La], [We], [E1]):
$$ 
	I(\pi)=\chi(S) + \chi_n(\pi), \quad 
	I_\pm(\pi)={1\over 2}\bigl( I(\pi)\pm c_1(X)\cdotp \pi(S)\bigr).  \eqno(1.1)
$$
The first of these formulas also holds if $S$ is unorientable.
Here $\chi_n(\pi)$ denotes the Euler number of the normal bundle 
$\nu_\pi =\pi^*(TX)/TS$ of the immersion $\pi$ (the self-intersection number 
of the zero section in $\nu_\pi$), and $c_1(X)\cdotp \pi(S)$ is the value of the 
first Chern class $c_1(X)= c_1(\Lambda^2 TX) \in H^2(X;Z)$ on the homology class 
$[\pi(S)]\in H_2(X;Z)$. The normal Euler number $\chi_n(\pi)$ is defined also 
for unorientable surfaces since any local orientation of $S$ coorients $\nu_\pi$ 
by the condition that the two orientations add up to the orientation of $X$ 
by the complex structure, and $\chi_n(\pi)$ is independent of these choices. 

Suppose now that $S$ is oriented and embedded in $X$. In this case its normal 
bundle $\nu$ is diffeomorphic to a \nbd\ of $S$ in $X$ and hence $\chi_n(S)$ 
is the {\it self-intersection number} $S^2$ of the homology class 
$[S]\in H_2(X;Z)$. Since $\chi(S)=2-2g(S)$, (1.1) gives
$$  
  I_\pm(S) = 1-g(S) + {1\over 2}\left( S^2 \pm c_1(X)\cdotp S\right)  \eqno(1.2)
$$
and the condition $I_\pm(S)\le 0$ in Theorem 1.1 is equivalent to 
$$
     g(S)\ge 1+ {1\over 2} \left( S^2+ |c_1(X)\cdotp S| \, \right).      \eqno(1.3)
$$ 
This is known as the {\it generalized adjunction inequality} and has been the 
subject of intensive recent research using the Seiberg-Witten approach to 
gauge theory. In the Appendix we summarize the results on (1.3) which are most 
interesting to us.  In particular, (1.3) holds if $X$ is a Stein surface and 
$S$ is not a homologically trivial two-sphere in $X$ (Theorem II in the Appendix).
Hence the condition $I_\pm (S)\le 0$ is also necessary for the existence 
of a Stein \nbd\ of $S$ in $X$. (It seems that this connection was first 
observed in [LM] and [N].) This gives

\proclaim 1.2 Corollary:   
Let $S$ be a closed, oriented, embedded real surface in a complex surface $X$. 
If $S$ is not a null-homologous two-sphere then $S$ is isotopic to an embedded
surface in $X$ with a regular Stein \nbd\ basis \iff\ (1.3) holds.

The adjunction inequality also holds in compact K\"ahler surfaces
with $b_2^+(X)>1$ (Theorem I in the Appendix). Hence Theorem 1.1
gives

\proclaim 1.3 Corollary: 
If $X$ is a compact K\"ahler surface with $b_2^+(X)>1$ then every closed, 
connected, oriented, embedded real surface $S\subset X$ of genus 
$g(S)>0$ can be $\cC^0$ approximated by an isotopic surface with a regular 
Stein \nbd\ basis. When $g(S)=0$ the same holds provided that $[S]\ne 0$
and none of the homology classes $\pm [S]$ can be represented by a 
(possibly reducible) complex curve in $X$.

Condition $b_2^+(X)>1$ in Corollary 1.3 cannot be omitted in general (Example 2 below). 
Consider now oriented surfaces $S\subset X$ without negative complex points. 
These include complex curves and, more generally, symplectic surfaces 
with respect to a positive symplectic form on $X$. 
(If $J\in {\rm End}(TX)$ denotes the almost complex structure on $X$
then a symplectic form $\omega$ on $X$ is $J$-{\it positive}, or $J$ is {\it tamed\/} 
by $\omega$, if $\omega(v,Jv)>0$ for any $0\ne v\in TX$. An immersion 
$\pi\colon S\to X$ is $\omega$-symplectic if $\pi^*\omega>0$ on $S$.) 
In this case (1.2) gives
$$ 
	0=2I_-(S)=\chi(S)+S^2-c_1(X)\cdotp S,
	\quad I_+(S)=c_1(X)\cdotp S = \chi(S)+S^2.  \eqno(1.4)
$$
Now $I_-(S)=0$ is equivalent to the  {\it genus formula} 
$g(S)=1+{1\over 2}(S^2-c_1(X)\cdotp S)$ (see [H, p.\ 361]),
and $I_+(S)\le 0$ \iff\ $c_1(X)\cdotp S\le 0$. 

\demo Example 1: If $S$ is a closed Riemann surface and $p\colon E\to S$ is a 
\holo\ line bundle then the zero section $S\subset E$ satisfies 
$I_+(S)= \chi(S)+ S^2 = 2-2g(S)+ c_1(E)$. (Here $c_1(E) \in Z$ is the 
value of the Chern class of $E$ on the fundamental class $[S]$ of $S$;
this is also called the {\it degree} of $E$.)
Hence the zero section can be approximated by an isotopic surface with 
a regular Stein \nbd\ basis in $E$ \iff\ $c_1(E)\le 2g(S)-2$.

%
%
\demo Example 2:  
Let $S$ be a closed, connected, oriented, embedded surface of degree 
$d\ge 1$ in $\CP^2$ (i.e., $[S]=d[H] \in H_2(\CP^2;Z)=Z$, where 
$H\simeq \CP^1$ is the projective line).  It is well known that 
$\L^2 (T\CP^2)=\cO_{\CP^2}(3)$ (its dual, the canonical bundle 
of $\CP^2$, is $K=\cO_{\CP^2}(-3)$). 
Hence $c_1(\CP^2)\cdotp [S]=c_1(\CP^2)\cdotp d[H]= 3d$
and (1.3) is equivalent to 
$$ 
	g(S)\ge 1+ {1\over 2}\left( S^2+ |c_1(\CP^2)\cdotp S| \right)= 
	 1+{1\over 2}(d^2+3d) ={1\over 2}(d+1)(d+2).         \eqno(1.5)
$$
In particular we must have $g(S)\ge 3$. 
If $S\subset\CP^2$ is isotopic to a complex or symplectic curve in $\CP^2$
then (by the genus formula) $g(S)={1\over 2}(d-1)(d-2)$ which is smaller than 
the right hand side in (1.5); hence such $S$ does not admit any Stein \nbd\ in $\CP^2$.

%
%
\demo Example 3:  Let $X=\CP^1\times\CP^1$.
For any $p\in \CP^1$ the group $H_2(X;Z)=Z\oplus Z$ is generated 
by $H_1=\CP^1\times \{p\}$ and $H_2=\{p\}\times \CP^1$. 
For an oriented embedded surface $S\subset X$
we let $d_j = S\cdotp H_j \in Z$ for $j=1,2$ denote the intersection
numbers with the two lines. The pair $d=(d_1,d_2)\in Z^2$ is called the 
{\it bidegree} of $S$. We have $S^2=2d_1d_2$, $c_1(X)\cdotp S =2(d_1+d_2)$,
and hence the generalized adjuntion inequality is
$$
	g(S) \ge 1+ d_1d_2 + |d_1+d_2|.                          \eqno(1.6)
$$
If (1.6) holds then by Theorem 1.1 $S$ is isotopic to a surface with a regular
Stein \nbd\ basis in $X$. Comparing (1.6) with the genus formula 
$g(C)=1+d_1d_2 - |d_1+d_2|$ for complex curves in $X$ (for which 
$d_1,d_2\ge 0$ and $d_1+d_2>0$) we see that none of them is isotopic
to a surface with a Stein \nbd. 
\endpr

Using (1.5) and elementary complex analysis S.\ Nemirovski proved that there 
exist no nonconstant \holo\ functions in any \nbd\ of such $S$ [N, Theorem 10], 
thus solving the {\it Vitushkin conjecture} to the effect that
there are no nonconstant holomorphic functions in any \nbd\ of 
an embedded sphere $S\subset \CP^2$ with $[S]\ne 0$. 
We show that (1.5) is sharp:

\proclaim 1.4 Proposition: For any $d\ge 1$ and $g\ge {1\over 2}(d+1)(d+2)$
there is an embedded oriented surface $S\subset \CP^2$ of genus $g$
and degree $d$ with a regular Stein \nbd\ basis.

For immersions in $\CP^2$ see Theorem 4.3 below. Proposition 1.4 is proved 
in section 3 and is a special case of the following.

\medskip\ni\bf  1.5 Theorem. \sl 
Let $S\subset X$ be a closed, connected, real surface smoothly embedded in 
a complex surface $X$. 
\item{(a)}  If $S$ is oriented and $k=\max\{I_+(S),I_-(S),0\}$ then 
$S$ is homologous to an embedded surface $S'\subset X$ of genus
$g(S')=g(S)+k$ with a regular Stein \nbd\ basis. 
\item{(b)} If $S$ is unorientable and  $k$ is the smallest integer with 
$3k\ge \max\{I(S),0\}$ then $S$ is $Z_2$-homologous in $X$ to an embedded 
unorientable surface $S'\subset X$ of genus $g(S')=g(S)+k$ 
with a regular Stein \nbd\ basis.
\medskip\rm

The next result shows that there are no genus restrictions for the existence 
of Stein \nbd s of embedded unorientable surfaces in $\CP^2$. 

\proclaim 1.6 Theorem: Every closed unorientable real surface embeds 
in $\CP^2$ with a regular basis of Stein \nbd s intersecting
every projective line.

\ni\it Remark. \rm By Theorem 1.8 in [F2] every closed unorientable surface 
$S$ also embeds in $\C^2$ with a regular Stein \nbd\ basis. The set of  
possible indices $I(\pi)$ of embeddings $\pi\colon S\hra\C^2$ 
equals $\{3\chi-4,3\chi,3\chi+4,\ldots,4-\chi\}$ where $\chi=2-g(S)$
[F2, p.\ 358]. This set contains both positive and negative numbers 
for any value of $\chi$. 
\endpr

Further results on Stein \nbd s of immersed surfaces are given 
in section 4. We mention a couple of related open problems.

%
%
\demo Problem 1:  Does there exist a Stein domain $\Omega\subset \C^2$ 
which is homotopically equivalent to the two-sphere~? This question, which
was raised in [F2], is apparently still open.  It is known that any embedded 
two-sphere $S\subset \Omega$ must be null-homologous in $\Omega$ 
[N, Theorem 15] and hence $[S]$ cannot generate $H_2(\CP^2;Z)$. 
There exists a precise description of the envelope of holomorphy of 
smooth two-spheres contained in closed \spsc\ hypersurfaces in $\C^2$ [BK, Kr].

\demo Problem 2: Let $S\subset X$ be an embedded real surface
with isolated complex points which admits a regular basis of Stein \nbd s. 
Is $I(p;S)\le 0$ for every $p\in S$? Clearly $S$ cannot have elliptic 
complex points due to Bishop discs [B, KW]. Further results on the existence 
of families of small analytic discs at certain complex points of higher index 
were obtained by Wiegerinck [Wg] and J\"oricke [J]. The problem seems 
open in general.

\beginsection 2. Construction of regular Stein \nbd s.

Theorem 2.1 below gives the existence of certain special immersions
regularly homotopic to a given initial immersion $S\to X$, and 
Theorem 2.2 provides regular Stein \nbd\ basis of special 
immersed surface.

\medskip \ni \bf 2.1 Theorem. \sl 
Any immersion $\pi_0\colon S\to X$ of a closed real surface $S$ into a complex 
surface $X$ can be $\cC^0$-approximated by a regularly homotopic smooth immersion 
$\pi\colon S\to X$ satisfying the following properties:
\item{(a)} At every complex point $p\in S$ of $\pi$ there are open \nbd s 
$p\in U\subset S$, $q=\pi(p)\in V \subset X$ and local \holo\ coordinates
$(z,w)$ on $V$ such that $z(q)=w(q)=0$ and $\pi(U) \subset V$ is given either 
by $w=z\bar z$ {\rm (a special elliptic point)} or by $w=\bar z^2$ 
{\rm (a special hyperbolic point)}.
If $I_\pm(\pi_0)\le 0$ (resp.\ $I(\pi_0)\le 0$ if $S$ is unorientable)
then we can choose $\pi$ as above to be without elliptic points. 
\item{(b)} The immersion $\pi$ only has transverse double points
(and no multiple points), and in a \nbd\ of each double point there exist 
local \holo\ coordinates $z=x+iy,\ w=u+iv$ on $X$ such that $\pi(S)$ is 
given locally by 
$ \{y=0,\ v=0\} \cup \{x=0,\ u=0\} = 
  (\R^2\times i\{0\}^2) \cup (\{0\}^2\times i\R^2)
$.
\item{(c)} Any embedding $S\hra X$ can be $\cC^0$-approximated by an isotopic 
embedding satisfying (a).
\medskip\rm

Theorem 2.1 is proved by cancellation of pairs of complex points due to 
Eliashberg and Harlamov (see the reference in [E2]). Their paper is not 
available in a standard source and we refer instead to Theorem 1.1 in [F2] 
where all details can be found. (We apologize for the incorrect reference 
[24] in [F2].) The result falls within the scope of Gromov's h-principle 
(sec.\ 2.4.5 in [Gr]). We recall the main steps.  

First we modify the immersion $\pi_0$ near each double point so that the new 
immersion satisfies property (b) (this is completely elementary and can be 
accomplished by a modification supported in small \nbd s of the double points). 
Suppose now that $p,q\in S$ are distinct isolated complex points of 
$\pi_0$ with $I(p;\pi_0)+I(q;\pi_0)=0$. If $S$ is oriented 
we assume that either both points are positive or both are negative. Choose a simple 
smooth arc $\gamma\subset S$ with endpoints  $p$ and $q$ which does not contain
any other complex point or double point of $\pi_0$. There exist \holo\ coordinates 
in an open set $U\supset \pi_0(\gamma)$ in $X$ which embed $U$ onto a domain in $\C^2$ 
such that $\pi_0(S)\cap U$ is mapped onto a graph 
$w=f(z)$ over a disc $D\subset \C$ (see section 5 in [F2] for the details). 
The winding number of ${\di\over \di \bar z}f$ around $\di D$ equals 
$I(p;\pi_0)+I(q;\pi_0)$ which is assumed to be zero. Such $f$ can be 
uniformly approximated by a smooth function $g$ on $D$ which 
equals $f$ near $\di D$ and satisfies ${\di\over \di\bar z}g\ne 0$ on $D$ 
(Lemma 4.1 in [F2] or Lemma 1 in [N, p.735]). This gives a $\cC^0$-approximation 
of $\pi_0$ by a regularly homotopic immersion which equals $\pi_0$ outside a small 
\nbd\ $V\supset \gamma$ in $S$ and is \tr\ on $V$. 

Using this procedure repeatedly one obtains an immersion $\pi\colon S\to X$ 
with $I_\pm(\pi)=I_\pm(\pi_0)$ such that each orientation class of $S$ contains 
only elliptic or only hyperbolic points of $\pi$, depending on the sign 
of $I_\pm(\pi)$. If $I_\pm(\pi)\le 0$ then all complex points of $\pi$ are hyperbolic. 
It is completely elementary to modify each elliptic resp.\ hyperbolic point to 
a complex point of special type as in Theorem 2.1 (a). 
If $\pi_0$ is an embedding, these deformations are carried out by 
an isotopy of embeddings. 
\endpr

%
%
\proclaim 2.2 Theorem: If $\pi\colon S\to X$ is an immersion satisfying Theorem 2.1 
and containing only special hyperbolic complex points then there are an open set $\Omega \subset X$ 
containing $M=\pi(S)$ and a smooth \psh\ function $\rho\colon \Omega\to [0,1]$ 
satisfying $M=\{x \in \Omega\colon  \rho(x)=0\}$, $d\rho\ne 0$ on $\Omega\bs M$, 
and such that for every $\e \in (0,1)$ the set 
$\Omega_\e =\{x\in \Omega\colon \rho(x)<\e \}$ is a smoothly 
bounded Stein domain which admits a strong deformation 
retraction onto $M$.

\demo Remarks: 1. An embedded real surface in a complex surface is locally 
holomorphically convex in a \nbd\ of any hyperbolic complex point [FSt]. 
From this it is easy to show that a surface $S\subset X$ with only hyperbolic 
complex points has a basis of Stein \nbd s (see [F2] and [N, Theorem 4]). 
However, at the time of this writing, {\it it is not known whether 
every such $S$ admits small Stein \nbd s homotopically equivalent to $S$}.
The problem is to find a `good' \psh\ function near every hyperbolic complex 
point whose critical points do not accumulate on $S$. (One must deal 
with a similar problem near each double point.) This is why we restrict our 
attention to {\it special hyperbolic complex points} and {\it special double points}. 

\ni 2. After the completion of this paper  M.\ Slapar informed 
the author that he had recently obtained regular Stein \nbd s 
in the presence of arbitrary hyperbolic complex points 
(personal communication, August 2002).

\demo Proof of Theorem 2.2: 
We first define $\rho$ near complex and double points. 
Let $p_1,\ldots,p_m\in M=\pi(S) \subset X$ be the complex points 
and $q_1,\ldots,q_k \in M$ the double points of $\pi$. By hypothesis there exist 
a \nbd\ $V_j\subset X$ of $p_j$ and \holo\ coordinates 
$\phi_j(p)=(z(p),w(p))$ on $V_j$ such that $z(p_j)=w(p_j)=0$,
$\phi_j(V_j)=r_j \disc\times \disc \subset \C^2$ for some $r_j\in(0,1)$
(where $\disc=\{\z\in \C\colon |\z|<1\}$), and $M\cap V_j =\{w=\bar z^2\}$. 
The nonnegative function  
$$ 
	\rho(z,w)=|w-\bar z^2|^2 = |w|^2+|z|^4-2\Re(wz^2)      
$$
has all the required properties in $V_j$. Indeed it is \psh\ (since 
$\Re(wz^2)$ is pluriharmonic), \spsh\ outside 
the complex disc $\Lambda_j=\{p\in V_j\colon z(p)=0\}$,
and $d\rho$ vanishes precisely on $M\cap V_j=\{w=\bar z^2\}$.

Each double point $q_j \in M$ of $\pi$ has a \nbd\ $W_j\subset X$ 
and \holo\ coordinates $\psi_j=(z,w)=(x+iy,u+iv)$ on $W_j$ such that 
in these coordinates $q_j=0$ and 
$M\cap W_j = \{y=0,\ v=0\} \cup \{x=0,\ u=0\}$. Set  
$$ 
    \rho(x+iy,u+iv)= (x^2+u^2)(y^2+v^2). 
$$
Its differential   
$$ 
	d\rho = 2 \bigl( x(y^2+v^2), y(x^2+u^2), u(y^2+v^2), v(x^2+u^2) \bigr) 
$$ 
vanishes precisely on $\{\rho=0\}=M\cap W_j$. We have 
$$  
	\rho_{z \bar z} = \rho_{w\bar w} = {1\over 2}(x^2+y^2+u^2+v^2), 
        \quad    \rho_{z\bar w} =  i(xv-yu) =\bar{\rho_{w\bar z}} \,.
$$ 
Since $\rho_{z \bar z}>0$ except at the origin, its complex Hessian $H_\rho(z,w)$ 
has at least one positive eigenvalue there. By Cauchy-Schwarz we have 
$$ 
   \eqalign{ 4\det H_\rho(z,w) &= (x^2+y^2+u^2+v^2)^2- 4(xv-yu)^2 \cr
   &\ge (x^2+y^2+u^2+v^2)^2 - 4(x^2+y^2)(u^2+v^2)  \cr
   &= |z|^4+|w|^4-2|z|^2|w|^2 \cr
   &= \bigl(|z|^2-|w|^2\bigr)^2. \cr}
$$
Thus $\det H_\rho(z,w)\ge 0$ which shows that both eigenvalues are nonnegative
and hence $\rho$ is \psh. The equality $\det H_\rho(z,w)=0$ holds precisely when 
$(x,y)=\l (v,-u)$ for some $\l\in \R$ and $|z|^2=|w|^2$, and this gives 
$w=\pm iz$. Let $L_j=\{p\in W_j\colon w(p)=\pm i z(p)\}$. Thus $\rho$ satisfies 
all required properties on $W_j$ and is \spsh\ on $W_j\bs L_j$. 

A function $\rho$ with the required properties has been defined on \nbd s 
$p_j\in V_j$ and $q_j\in W_j$. We now extend $\rho$ to a smooth nonnegative function 
in a \nbd\ $V\supset M$ such that the extension vanishes precisely on $M$ 
and its real Hessian is nondegenerate in any normal direction to $M$ 
at all points of  $M_0=M\bs \{p_j,q_j\}$. More precisely, if 
we denote by $\nu\to M_0$ the normal bundle to $M_0$ in $X$ and realize it as a 
subbundle of $TX|_{M_0}$ such that $TX|_{M_0}=T M_0\oplus \nu$, we require
that the second order derivatives of $\rho$ in the fiber directions
$\nu_x$ at any $x\in M_0$ give a nondegenerate quadratic form on $\nu_x$
(which is necessarily positive definite since $\rho$ has a local minimum at 
$0_x\in \nu_x$). One can obtain such an extension by taking a suitable 
second order jet along $M$ with the required properties and applying 
Whitney's theorem to find a function which matches this jet. 

A more explicit construction is the following. Choose a Riemannian metric 
$h(x;\xi,\eta)$ on the normal bundle $\nu=\{(x,\xi)\colon x\in M_0,\xi\in \nu_x\}$. 
The function $\wt\rho(x,\xi) = h(x;\xi,\xi)\in \R_+$ has suitable properties on $\nu$ 
(with $M_0$ corresponding to the zero section of $\nu$). Let $\psi$ be a 
diffeomorphism of $\nu$ onto a tubular \nbd\ $U\subset X$ of $M_0$ such 
that $\psi(x,0)=x$ and $d\psi(x,0)$ is the identity for each $x\in M_0$. 
Then $\rho_0=\psi\circ\wt\rho \circ \psi^{-1}\colon U\to \R_+$
satisfies the required properties near $M_0$.
Within each coordinate \nbd\ $V_j$ and $W_j$ as above we patch 
$\rho_0$ with the previously chosen function $\rho=\rho_j$ on this 
\nbd\ by taking $\rho=\chi_j \rho_j+ (1-\chi_j)\rho_0$ on $V_j$,
where $\chi_j$ is a smooth cut-off function which is supported in $V_j$ 
(resp.\ in $W_j$) and equals one in a smaller \nbd\ of $p_j$ resp.\ $q_j$. 
At points $x\in M_0 \cap V_j$ the real Hessian of $\rho$ equals   
$$
    H^\R_\rho(x) = \chi_j(x) H^\R_{\rho_j}(x) +(1-\chi_j(x)) H^\R_{\rho_0}(x). 
$$
This is nonnegative on $T_x X$ and strongly positive on $\nu_x$ 
(since it is a convex combination of two forms with these properties). 

We claim that the function $\rho$ obtained in this way is \psh\ in 
a \nbd\ of $M$, it has no critical points near $M$ 
(except on $M$), and the sublevel sets $\{\rho<\e\}$ 
for sufficiently small $\e>0$ are Stein. We first show 
plurisubharmonicity. By construction $\rho$ is such near 
each $p_j$ and $q_j$. Let $p\in M_0 =M\bs \{p_j,q_j\}$. 
Since $M_0$ is \tr, there exist local \holo\ coordinates 
$(z,w)=(x+iy,u+iv)$ near $p$ such that in these coordinates $p=0$ and 
$T_p M_0=\R^2\oplus i\{0\}^2 =\{y=0,\ v=0\}$. Since $\rho=0$ and 
$d\rho=0$ on $M_0$, those second order derivatives of $\rho$ at $0$ 
which contain at least one differentiation on $x$ or $u$ 
vanish at $0$ and hence
$$ 4\rho_{z\bar z}(0)= \rho_{yy}(0), \quad 
   4\rho_{w\bar w}(0)= \rho_{vv}(0), \quad
   4\rho_{z\bar w}(0)= \rho_{yv}(0).
$$
Thus  the complex Hessian of $\rho$ at $0$ equals one quarter of 
the real Hessian of the function $(y,v)\to \rho(0+iy,0+iv)$ at $y=0,v=0$.
By construction this is positive definite and hence
$\rho$ is \spsh\ at every $p\in M_0$. Thus $\rho$ is \psh\ 
in a small \nbd\ $\Omega\supset M$ and \spsh\ in 
$\Omega'=\Omega\bs \bigl((\cup \L_j) \cup (\cup L_j) \bigr)$.
We may choose $\Omega=\{\rho <\e_0\}$ for some $\e_0>0$. 

Next we show that $d\rho \ne 0$ on $\Omega\bs M$ if $\e_0>0$ (and hence $\Omega$) 
are chosen sufficiently small. By construction this holds in small \nbd s 
of the points $p_j$ and $q_j$. Over $M_0=M\bs \{p_j,q_j\}$
we consider the conjugate function $\wt\rho=\psi^{-1}\circ \rho\circ \psi$
on the normal bundle $\nu\to M_0$. Suppose that $d\wt\rho(x,\xi)\cdotp\xi=0$
for some $(x,\xi)\in \nu$ with $x\in M_0$ and $\xi\ne 0$. Consider
the function $t\in\R\to \wt\rho(x,t\xi)$. By hypothesis its derivative
vanishes at $t=0$ and $t=1$ and hence its second derivative vanishes at 
some $t_0 \in (0,1)$. This means that the Hessian of $\wt\rho(x,\cdotp)$ 
vanishes at the fiber point $t_0\xi \in \nu_x$ in the direction of the 
vector $\xi \in \nu_x$. We have seen that this does not happen in a sufficiently 
small \nbd\ of the zero section of $\nu$ which establishes our claim.

Thus for any $\e\in (0,\e_0)$ the set 
$\Omega_\e=\{ x\in \Omega \colon \rho(x)<\e\}$ is a relatively
compact smoothly bounded (weakly) pseudoconvex domain in $X$.
We obtain a strong deformation retraction of $\Omega_\e$
onto $M$ by integrating the flow of the negative gradient
of $\rho$ from $t=0$ to $t=+\infty$ and rescaling the time
interval to $[0,1]$. (The gradient of $\rho$ is defined by
the equation $\nabla \rho \rfloor h=d\rho$ where $h$ is a
Riemannian metric on $X$ and $\rfloor$ denotes the contraction.) 

It remains to show that $\Omega_\e$ is Stein for each $\e\in (0,\e_0)$. 
Let $h_\e\colon (-\infty,\e) \to \R$ be an increasing strongly 
convex function with $\lim_{t\to \e}h(t)=+\infty$.
If there exists a \spsh\ function $\psi$ in a \nbd\ of 
$\bar\Omega_\e$ in $X$ (this is the case for instance if $X$
is Stein) then $\psi + h_\e\circ \rho$ is a \spsh\ exhaustion 
function on $\Omega_\e$ and hence $\Omega_\e$ is Stein according
to [Gra]. 

In general a weakly \psc\ domain in a non-Stein manifold 
need not be Stein. In our situation we proceed as follows. 
Recall that the Levi form of $\rho$ and hence of $\rho_\e=h_\e\circ \rho$ 
is degenerate only on the complex curves $\L_j \subset \Omega$ and 
$L_j\subset \Omega$. These curves intersect $M$ transversely at the points 
$p_j$ resp.\ $q_j$. We take $\phi=\rho_\e +\d\tau$ where $\d>0$ 
and $\tau$ is a smooth function in a \nbd\ of $\bar \Omega_\e$ 
which is \spsh\ on the curves $\L_j$ and $L_j$ and which vanishes 
outside the coordinate \nbd\ $V_j$ of $p_j$ (resp.\ $q_j$). Since 
the complex Hessian of $\rho_\e$ is bounded away from zero on the 
set in $\Omega_\e$ where $\tau$ fails to be \psh, $\phi$ is a 
\spsh\ exhaustion function on $\Omega_\e$ provided that $\d>0$ is 
chosen sufficiently small. 

We take $\tau$ to be given in local coordinates $(z,w)$ near
$p_j$ resp.\ $q_j$ by $\tau(z,w)=\chi(z,w)(|z|^2+|w|^2)$
where $\chi$ is a suitable cut-off function. At $p_j$
we have $\L_j=\{z=0\}$, the coordinate \nbd\ $V_j\subset X$ of $p_j$ 
is mapped onto $\{|z|<r,\ |w|<1\} \subset \C^2$, and a suitable 
cut-off function is $\chi(|z|)$ where $\chi(t)=1$ for $t\le r/2$ and 
$\chi(t)=0$ for $t\ge 3r/4$. At $q_j$ we may assume that the coordinate
\nbd\ $W_j\subset X$ is mapped onto $\{|z|<1,\ |w|<1\} \subset \C^2$
and in this case a suitable cut-off function is $\chi(|z|^2+|w|^2)$
where $\chi(t)=1$ for $t\le 1/4$ and $\chi(t)=0$ for $t\ge 3/4$.
In both cases the support of the differential $d\chi$ intersects $\Omega_\e$
for all sufficiently small $\e>0$ in a set whose closure (in $X$)
is compact and does not meet any of the curves $\L_j$ and $L_j$. 
On this set the eigenvalues of the Levi form of $\rho_\e=h_\e\circ \rho$ 
are bounded away from zero and hence for $\d>0$ sufficiently small 
$\rho_\e+\d\tau$ is a \spsh\ exhaustion function in $\Omega_\e$. 
Hence $\Omega_\e$ is Stein by Grauert's theorem [Gra].

Note that we can perturb each domain $\Omega_\e$ constructed above 
to a \spsc\ domain in $X$ by adding to $\rho$ a small function which 
is \spsh\ on the sets $\L_j\cap b\Omega_\e$ and $L_j\cap b\Omega_\e$. 
\endpr

\ni \it Proof of Theorem 1.1. \rm
Let $S\subset X$ be an embedded surface 
with $I_\pm(S)\le 0$ (resp.\ $I(S)\le 0$ if $S$ is unorientable).
Let $S'\subset X$ be an isotopic embedding satisfying the conclusion
of Theorem 2.1. Theorem 2.2 then gives a regular basis of Stein \nbd s 
of $S'$. The analogous conclusions hold for immersions.

%
%
%
%
\beginsection 3. Connected sums of embedded and immersed surfaces.

Let $X_1$ and $X_2$ be smooth $n$-manifolds. Choose embedded $n$-dimensional
discs $D_j\subset X_j$ for $j=1,2$ and let $\phi\colon \di D_1\to \di D_2$
be a smooth diffeomorphism which is orientation reversing if the two
manifolds are oriented. The identification space 
$(X_1\bs {\rm int}D_1)\cup_\phi (X_2\bs {\rm int}D_2)$ can be given the 
structure of a smooth manifold called the 
{\it connected sum} $X_1\# X_2$ of $X_1$ and $X_2$ (see [GS], p.\ 20). 
The smoothing process can be visualized by connecting the two pieces 
$X_j\bs {\rm int} D_j$ by a cylinder $[0,1]\times S^{n-1}$  glued 
along its boundary spheres $(\{0\}\times S^{n-1}) \cup (\{1\}\times S^{n-1})$ 
to the spheres $\di D_1$ resp.\ $\di D_2$. If $X_1=X_2=X$ and 
we wish to use as $\phi$ the identity map on $\di D \subset X$, 
we must reverse  the orientation on one of the copies of $X$; 
in this case we write $X\# \bar X$. For example, if $X$ is a complex 
$n$-manifold then $X\# \bar{\CP^n}$ is diffeomorphic to the complex 
manifold obtained by blowing up a point in $X$. The connected sum operation 
extends to several terms and we shall write $X_1\# kX_2$ for the connected sum 
of $X_1$ with $k$ copies of $X_2$.

Suppose now that $S_1,\ S_2\subset X$ are embedded or immersed real surfaces 
in general position in a complex surface $X$. 
We can realize their connected sum $S_1\# S_2$ as an immersed 
surface in $X$ with the same number of double points and with 
two additional hyperbolic complex points, one positive and one negative
if the surfaces are oriented, which are contained in the attaching
cylinder $\Sigma \simeq [0,1]\times S^1$ (see Sect.\ 3 in [F2] for the details). 
The two removed discs $D_j \subset S_j$ are chosen  totally real. 
Thus we have $I(S_1\# S_2)=I(S_1)+I(S_2)-2$, and if both surfaces are oriented 
we also have 
$$ 
	I_+(S_1\# S_2)=I_+(S_1)+ I_+(S_2)-1,\quad  
	I_-(S_1\# S_2)=I_-(S_1)+ I_-(S_2)-1.
$$
Furthermore, $g(S_1\# S_2)= g(S_1)+g(S_2)$ and 
$[S_1\# S_2] = [S_1]+[S_2] \in H_2(X;Z)$. 
If $S_1$ and $S_2$ are disjoint embedded surfaces in $X$ then 
$S_1\# S_2$ can also be realized as an embedded surface. 

\demo Proof of Theorem 1.5:
Let $T\subset X$ be an embedded null-homologous \tr\ torus
(take a small torus in a coordinate chart). 
For any embedded oriented surface $S\subset X$ we have
$$
   g(S\#T)=g(S)+1,\quad I_\pm(S\#T)=I_\pm(S)-1,\quad [S]=[S\#T]\in H_2(X;Z).
$$
Attaching $k=\max\{I_+(S),I_-(S),0\}$ torus handles we obtain an embedded surface 
$S\# kT \subset X$ with $g(S\#kT)=g(S)+k$ and $I_\pm(S\#kT)=I_\pm(S)-k \le 0$. 
Part (a) now follows from Theorem 1.1. The same argument applies to immersed
surfaces.

Suppose now that $S\subset X$ is an embedded (or immersed) unorientable surface. 
Let $M\subset X$ be a copy of the real projective plane $\RP^2$ 
embedded in a coordinate chart in $X$ with $I(M)=-1$ [F2, p.\ 367]. 
Then $g(S\#M)=g(S)+1$ and $I(S\# M)=I(S)-3$. 
Thus $g(S\#kM)=g(S)+k$, $I(S\# kM)=I(S)-3k\le 0$, and hence (b) follows 
from Theorem 1.1. We can use other types of handles such 
as a \tr\ torus $T$ or a Klein bottle $K$ (for an explicit \tr\ Klein bottle 
in $\C^2$ see [R]). The surfaces $S\# T$ and $S\# K$ are diffeomorphic and 
satisfy $I(S\# T)=I(S)-2$ and $g(S\# T)=g(S)+2$; hence these handles are 
less effective in lowering the index than $\RP^2$.
\endpr

\demo Proof of Proposition 1.4:  
Let $C\subset \CP^2$ be a smooth complex curve of degree $d$.
By (1.4) we have $I_-(C)=0$, $I_+(C)=3d$, $g(C)=g_\C(d)={1\over 2}(d-1)(d-2)$.
Attaching $k \ge 3d$ torus handles we obtain an embedded surface 
$S=C\#k T \subset \CP^2$ homologous to $C$, with $I_\pm(S)\le 0$ and 
$$
   g(S)=g(C)+k = {1\over 2}(d-1)(d-2)+k \ge {1\over 2}(d+1)(d+2).   
$$
The result now follows from Theorem 1.1. Note that, by the Thom conjecture 
proved in [KM], a smooth complex curve of degree $d$ in $\CP^2$ has the smallest 
genus $g_\C(d)={1\over 2}(d-1)(d-2)$ among all smooth oriented real surfaces 
of degree $d$ in $\CP^2$. For the symplectic case see [OS] and [MST].
\endpr

\demo Proof of Theorem 1.6:  Let $\RP^2 \simeq M\hra \C^2 \subset\CP^2$ be 
an embedding with index $I(M)=-1$ mentioned above [F2, p.\ 367]. 
Let $S=\{[x\colon y\colon z]\in\CP^2 \colon x,y,z\in\R\}\simeq\RP^2$.
Clearly $S$ is \tr\ and intersects every projective line in $\CP^2$;
in fact their (mod 2) intersection number  equals one.
The connected sum $S\# kM \subset \CP^2$ is an embedded unoriented surface 
with genus $1+k$ and index $-3k\le 0$, and hence Theorem 1.1 applies.

%
%
\beginsection 4. Regular Stein \nbd s of immersed surfaces.

In this section $S$ denotes a closed connected oriented real surface 
of genus $g(S)$. We shall consider immersions $\pi\colon S\to X$ into a complex 
surface $X$ with simple (transverse) double points and with no multiple 
points.  Furthermore we shall assume that both tangent planes at any double 
point are \tr. At each double point $\pi$ has  self-intersection index 
$\pm 1$ which is independent of the choice of the orientation on $S$. 
Double points with index $+1$ will be called {\it positive}
and those with index $-1$ will be called {\it negative}. If $\pi$ has 
$\d_+$ positive and $\d_-$ negative double points then
$\d(\pi)=\d_+ - \d_-$ is the (geometrical) {\it self-intersection index} 
of $\pi$ which only depends on its regular homotopy class.

Recall that $\chi_n(\pi)$ denotes the normal Euler number of $\pi$. It is easily 
seen that each double point of $\pi$ contributes $\pm 2$ (the sign depending on 
its self-intersection index) to the homological self-intersection number 
of the image $\pi(S)$ in $X$. This gives $\chi_n(\pi)+2\d(\pi)=\pi(S)^2$. 
From (1.1) we get  
$$ 
	I_\pm(\pi)=1-g(S)- \d(\pi) + 
	{1\over 2}\left( \pi(S)^2 \pm  c_1(X)\cdotp \pi(S)\right), 
							                \eqno(4.1)
$$
and the condition $I_\pm(\pi)\le 0$ is equivalent to 
$$ 
	g(S)+ \d_+ 
	\ge 1 + \d_- + {1\over 2}(\pi(S)^2+|c_1(X)\cdotp \pi(S)|\, ).      \eqno(4.2) 
$$

\proclaim 4.1 Corollary: If (4.2) holds then $\pi$ is regularly homotopic to an 
immersion $S\to X$ whose image has a regular Stein \nbd\ basis in $X$.
The regular homotopy can be chosen such that it preserves the location
(and number) of double points. Conversely, if $\pi(S)$ has a Stein 
\nbd\ $\Omega\subset X$ and if $[\pi(S)]\ne 0$ in $H_2(\Omega;Z)$ then
the following generalized adjunction inequality holds:
$$ 
   g(S) + \d_+ \ge 1 +{1\over 2}(\pi(S)^2+|c_1(X)\cdotp \pi(S)|).       \eqno(4.3) 
$$

Unlike in (4.2), the term $\d_-$ is not present on the right hand side of (4.3)
and hence the two inequalities coincide only when $\pi$ has no negative double points.
If $\d_->0$, there is a `gray area' between the direct and the converse part
of Corollary 4.1. The direct part follows from Theorems 2.1 and 2.2; 
for the converse part see Theorem III in the Appendix. 
When $\pi$ has no negative complex points (for instance 
if $\pi$ is complex or symplectic) then $I_-(\pi)=0$ which gives 
the {\it genus formula} (see e.g.\ [IS, p.\ 576])
$$ 
	g(S)=1-\d(\pi) + {1\over 2}\left( \pi(S)^2 - c_1(X)\cdotp \pi(S)\right).
$$

\proclaim 4.2 Theorem: Assume that $S$ is an oriented closed surface  
and $\pi \colon S\to X$ is an immersion with $k$ double points.  If 
$$  
	g(S)+ k \ge 1 +{1\over 2}(\pi(S)^2+|c_1(X)\cdotp \pi(S)|\, )     \eqno(4.4)
$$
then every open \nbd\ $\Omega\subset X$ of $\pi(S)$ contains an embedded 
oriented surface $M\subset \Omega$ of genus $g(M)=g(S)+ k$
such that $[M]=[\pi(S)] \in H_2(\Omega;Z)$ and $M$ admits a 
regular basis of Stein \nbd s.

\demo Proof: 
At each double point of $\pi(S)$ we replace a pair of small intersecting 
discs in $\pi(S)$ by an embedded annulus $\Sigma \simeq S^1\times[0,1]$. 
This well known procedure, which amounts to replacing double points by 
handles, can be done within $\Omega$ and it does not change the homology 
class of the image, but it increases the genus of the immersed surface by $k$. 
Here is a precise description.

By Theorem 2.1 we may assume that in local \holo\ coordinates 
$(z,w)=(x+iy,u+iv)$ on $X$ near a double point the immersed surface 
is the union of discs in the lagrangian planes 
$\Lambda_1 = \{y=0,\ v=0\} \subset \C^2$, $\Lambda_2 =\{x=0,\ u=0\} \subset\C^2$. 
Let $\Lambda_1$ be oriented by ${\di_x}\wedge {\di_u}$ and 
$\Lambda_2$ by $\kappa\, {\di_v}\wedge {\di_y}$ where $\kappa=\pm 1$ 
is the self-intersection index of the double point.  

If $\kappa =+1$ we can use the handle $\Sigma_+=\{(x+iu)(y-iv)=\e\} \cap B$ 
for a small $\e\ne 0$, where $B \subset\C^2$ is a small ball around the origin. 
Outside of $B$ we smoothly patch $\Sigma_+$ with 
$(\Lambda_1\bs D_1)\cup (\Lambda_2 \bs D_2)$ without introducing new complex points 
($D_j$ is a disc in $\L_j$). A simple calculation shows that $\Sigma_+$ 
is \tr\ in $\C^2$ for every $\e\ne 0$, and hence this handle does 
not change the indices $I_\pm$. 

If $\kappa= -1$, an appropriate handle is $\Sigma_-=\{(x+iu)(y+iv)=\e\}$ 
for small $\e\ne 0$. It has four hyperbolic complex points at 
$x=y =\pm \sqrt{\e/ 2}$, $v=-u= \pm \sqrt{\e/ 2}$ (independent choices of signs),
two positive and two negative. Hence $I_\pm(\Sigma_-)=-2$ and the indices 
of the immersion decrease by two. 

After replacing all double points by handles in this way we obtain 
an embedded surface $M\subset \Omega$ with $g(M)=g(S)+k$
and $[M]=[\pi(S)]\in H_2(X;Z)$. From (4.4) it follows that 
$g(M) \ge 1 +{1\over 2}(M^2+|c_1(X)\cdotp M|\, )$ and hence 
Theorem 1.1 applies to $M$.
\endpr

\demo Remark: The replacement of double points by handles can also 
be explained using the following complex model. 
In local \holo\ coordinates $(z,w)=(x+iy,u+iv)$ we represent 
the double point by $\{zw=0\}= L_1\cup L_2$, where $L_1= \{w=0\}$ is 
oriented by $\di_x\wedge \di_y$ and $L_2=\{z=0\}$ is oriented by 
$\kappa\, \di_u\wedge \di_v$ where $\kappa=\pm 1$ is the 
intersection index. An appropriate handle is $zw=\e$
when $\kappa=+1$ and $z\bar w=\e$ when $\kappa=-1$.
\endpr

We now consider immersed surfaces in $\CP^2$. Recall that the degree of  
$\pi\colon S\to \CP^2$ is the integer $d$ satisfying 
$[\pi(S)]=d[H] \in H_2(\CP^2;Z)=Z$ where $H$ is the projective 
line. Reversing the orientation on $S$ if necessary we may assume $d\ge 0$. 
The number on the right hand side of (4.2)--(4.4) equals 
$1+{1\over 2}(d^2+3d)={1\over 2}(d+1)(d+2)$.

%
%
\proclaim 4.3 Theorem: 
Let $S_g$ be a closed oriented surface of genus $g$ and let 
$\pi\colon S_g\to \CP^2$ be an immersion of degree $d\ge 1$ 
with $\d_+$ positive double points. If $\pi(S_g)$ has a Stein 
\nbd\ in $\CP^2$ then $g+\d_+ \ge {1\over 2}(d+1)(d+2)$. Conversely, for any pair 
of integers $g, \d_+ \ge 0$ satisfying the above inequality there exists an immersion
$\pi\colon S_g \to \CP^2$ of degree $d$ with $\d_+$ positive double points 
(and no negative double points) whose image $\pi(S_g)$ admits a regular Stein 
\nbd\ basis in $\CP^2$. In particular, there is an immersed sphere in $\CP^2$
in the homology class of the projective line, with three double points, 
which admits a regular Stein \nbd\ basis.

\demo Proof: The first part follows from Corollary 4.1. We prove the 
converse part by an explicit construction. 
Let $\pi_0\colon S_0\to \C^2$ be Weinstein's  lagrangian (hence \tr) 
`figure eight' immersion of the sphere [Wn]:
$$ \pi_0(x,y,u)=\bigl( x(1+2iu),y(1+2iu)\bigr),\quad x^2+y^2+u^2=1. $$
It identifies the points $(0,0,\pm1)$ and the corresponding double point 
at $0\in\C^2$ has self-intersection index $+1$.

Let $C_1,C_2,\ldots,C_d\subset \CP^2$ be projective lines in general position. 
There are ${1\over 2}d(d-1)$ positive intersection points. For each 
$j\in\{2,\ldots,d\}$ we replace one of the intersection points of $C_j$ with 
$C_1\cup\cdots\cup C_{j-1}$ by a handle as in the proof of Theorem 4.2. 
This eliminates $d-1$ double points and changes the given system of lines to 
an immersed sphere $\rho\colon S_0\to\CP^2$ of degree $d$ with $I_+(\rho)=3d$, 
$I_-(\rho)=0$ and with ${1\over 2}(d-1)(d-2)$ positive double points. 

Write $g+\d_+=l+{1\over 2}(d+1)(d+2)$, with $l\ge 0$, and set $k=l+3d\ge 3d$.
Take $k$ pairwise disjoint copies $M_1,\ldots,M_k \subset \CP^2\bs \rho(S_0)$ 
of Weinstein's sphere $\pi_0(S_0)$. The connected sum 
$\rho(S_0)\# M_1\#\cdots\# M_k \subset \CP^2$ 
is parametrized by a degree $d$ immersion $\pi\colon S_0\to \CP^2$  with 
$I_+(\pi)=3d-k\le 0$, $I_-(\pi)=-k< 0$, and with $k+{1\over 2}(d-1)(d-2) = g+ \d_+$ 
positive double points. Replacing $g$ of the double points by handles we obtain an 
immersed surface of genus $g$ with $\d_+$ positive double points and 
with non-positive indices $I_\pm$. The result now follows from Theorem 1.1.
\endpr

By a theorem of Gromov [Gr, p.\ 334, Theorem A] every immersion 
$\pi\colon S\to \CP^2$ of degree $d\ge 1$ can be $\cC^0$-approximated
by a symplectic immersion $\wt \pi\colon S\to \CP^2$, i.e., 
$\wt\pi^*(\omega_{FS})>0$ on $S$ where $\omega_{FS}$ denotes the
Fubini-Studi form on $\CP^2$. (If we specify a two-form $\theta>0$ 
on $S$ with $\int_{S}\theta =1$, we can even choose $\wt\pi$ such that 
$\wt\pi^* \omega_{FS}= d\cdotp \theta$.) Hence Theorem 4.3 implies

\proclaim 4.4 Corollary: For any $d\ge 1$ there exists a symplectically 
immersed sphere in $\CP^2$ of degree $d$ with a Stein \nbd.

I thank S.\ Ivashkovich for telling me the question answered in part by 
Corollary 4.4 and for pointing out the relevance of Gromov's theorem. 
The results of [IS] imply that {\it there exist no immersed symplectic spheres 
in $\CP^2$ with a Stein \nbd\ basis}.

Corollary 4.4 may seem somewhat surprising since symplectic curves share 
many properties with complex curves. The main difference (which is 
essential here) is that complex curves only intersect positively 
while symplectic curves may also intersect negatively.

The symplectic approximation furnished by Gromov's theorem need not 
be regularly homotopic to the initial immersion and hence the number of double
points may increase. The question remains about the minimal number of 
double points of immersed symplectic spheres in $\CP^2$ of a given degree 
$d\ge1$ which admit Stein \nbd s. For $d=1$ the genus formula gives 
$\d(\pi)=1+{1\over 2}(d^2 - 3d)=0$ which shows that the number of positive
double points equals the number of negative double points. 
By Theorem 4.3 we have $\d_+ \ge {1\over 2}(d+1)(d+2)=3$ 
and hence {\it there are at least $6$ double points}.

\demo Question: Is there an immersed symplectic sphere in $\CP^2$ 
of degree one (in the homology class of the projective line) with six double 
points and with a Stein \nbd~?

\vfill\eject

%
%
%
%
%
\medskip\ni\bf 
Appendix: Generalized adjunction inequalities.
\medskip\rm

We say that a closed, connected, oriented, smoothly embedded real surface 
$S$ in a complex surface $X$ satisfies the 
{\it generalized adjunction inequality} if
$$  
	g(S)\ge 1+ {1\over 2}\left( S^2+|c_1(X)\cdotp S| \right),     \eqno(*) 
$$
where $g(S)$ is the genus of $S$, $S^2$ is the self-intersection number
of its oriented homology class $[S]\in H_2(X;Z)$ in $X$, 
and $c_1(X)=c_1(\L^2 TX) \in H^2(X;Z)$ is the first Chern class 
of $X$. The purpose of this appendix is to collect from the literature 
those results on (*) which are most relevant to this paper. 
We emphasize that no originality whatsoever is being claimed on our part.

If $X$ is a closed oriented four-manifold then $H^2(X;\R)$ is equipped with 
the intersection pairing induced by $(\a,\b)\to \int_X \a\wedge \b$ where 
$\a,\b$ are closed 2-forms on $X$. This pairing is Poincar\'e dual to the  
intersection pairing of homology classes in $H_2(X;Z)$. Let $b_2^+(X)$ 
denote the dimension of a maximal linear subspace of $H^2(X;\R)$ on which 
this pairing is positive definite. (If $X$ is a compact K\"ahler surface then 
$b_2^+(X)>0$ is an odd positive integer.)

%
%
\medskip\ni \bf Theorem I. {\rm ([KM], [MST], [FiS], [OS])}  \sl
Assume that $X$ is a compact K\"ahler surface with $b_2^+(X)>1$.
Let $S\subset X$ be a smoothly embedded, closed, connected,
oriented real surface in $X$. Then the generalized adjunction 
inequality (*) holds in each of the following cases:
\item{(a)}   $g(S)>0$ (i.e., $S$ is not a sphere), 
\item{(b)}   $g(S)=0$, $[S]\ne 0$ in $H_2(X;Z)$, 
and none of the homology classes $\pm[S]$ can be represented 
by a (possibly reducible) complex curve in $X$.
\rm

\medskip
For homologically trivial surfaces (*) reduces to $g(S)\ge 1$ which only 
excludes the sphere. For homologically nontrivial spheres (*) requires 
in particular that $S^2\le -2$ which fails for the exceptional sphere 
of a blown-up point (this has $S^2=-1$). The inequality (*) fails for any complex 
(or symplectic) curve $C\subset X=\CP^2$ since $g(C)={1\over 2}(d-1)(d-2)$
by the classical genus formula (where $d>0$ is the degree of $C$), but the right 
hand side of (*) equals ${1\over 2}(d+1)(d+2)$. Hence the condition 
$b^+_2(X)>1$ cannot be omitted.

Before proceeding we need to recall the notion of the Seiberg-Witten function
(see [GS, Sect.\ 2.4] for a concise exposition and [M] for a more 
comprehensive one). If $X$ is a simply connected, closed, oriented four-manifold
with $b_2^+(X)>1$, the {\it Seiberg-Witten function} on $X$ is an integer 
valued function $SW_X\colon \cC_X\to Z$ defined on the set of all characteristic 
classes $\cC_X \subset H^2(X;Z)$ (see [GS, p.\ 51]). One of its main features is that 
the value $SW_X(K)$ on any $K\in \cC_X$ only depends on the smooth
structure on $X$. The story is more complicated if $X$ is not simply 
connected or if $b_2^+(X)=1$ and we shall not go into that. A characteristic class 
$K$ is a {\it Seiberg-Witten basic class} if $SW_X(K)\ne 0$, and $X$ is 
of {\it Seiberg-Witten simple type\/} if for any such $K$ the associated moduli 
space of solutions of the perturbed Seiberg-Witten equations is zero dimensional
(hence finite). 

For surfaces with $S^2\ge 0$ Theorem I was proved by Kronheimer and Mrowka [KM] 
and Morgan, Szab\'o and Taubes [MST] when $X$ is an oriented closed four-manifold 
with $b_2^+(X)>1$, of Seiberg-Witten simple type, and with $c_1(X)$ 
replaced in (*) by any Seiberg-Witten basic class.  
Since any compact symplectic four-manifold (hence any compact 
K\"ahler surface) with $b_2^+(X)>1$ is of Seiberg-Witten simple type [T2] 
and satisfies $SW_X(c_1(X))\ne 0$ [T1], Theorem I is a special case
of the quoted results. Note that $S$ cannot be a sphere since the existence 
of an embedded sphere with $S^2\ge 0$ and $[S]\ne 0$ implies that 
$SW_X$ vanishes identically [FiS].

In the case $S^2<0$ Theorem I was proved for spheres by Fintushel 
and Stern [FiS, Theorem 1.3] and for surfaces 
of positive genus by Ozsv\'ath and Szab\'o [OS, Corollary 1.7].
The proofs essentially depend on the assumption that $X$ is of  
Seiberg-Witten simple type, but the K\"ahler hypothesis can be relaxed 
as above. The argument goes as follows.  If (*) fails for a given $S\subset X$ 
then $K=c_1(X)+2\e PD(S) \in H^2(X;Z)$ is a Seiberg-Witten basic class;
here $\e=\pm 1$ is the sign of $c_1(X)\cdotp S$ and $PD(S)$ is the 
Poincar\'e dual of $S$ (Theorem 1.3 in [OS]; Remark 1.6 in [OS] explains 
the factor $2$ in front of $PD(S)$ as opposed to the notation in Theorem 1.3). 
Furthermore, if $g(S)>0$, the corresponding moduli space of solutions of the 
perturbed Seiberg-Witten equations has positive dimension (see the proof of 
Corollary 1.7 in [OS]) which contradicts the simple type assumption on $X$. 
If $g(S)=0$ then one cannot draw the last conclusion. 
However, $SW_X(K)\ne 0$ for a $K\in H^2(X;Z)$ implies that the class 
$K'={1\over 2}(c_1(X)-K)$ is Poincar\'e dual to a complex curve $C\subset X$ 
[GS, p.\ 417]. Applying this to $K=c_1(X)+2\e PD(S)$ for which $K'= \pm PD(S)$ 
we obtain $\pm[S]=[C]$ as claimed (see [FiS, Theorem 1.3]).
\endpr

\proclaim Theorem II: {\bf (Adjunction inequalities in Stein surfaces)} 
Let $X$ be a Stein surface and let $S\subset X$ be a smoothly embedded, closed,
connected, oriented real surface in $X$. Then the generalized adjunction 
inequality (*) holds unless $S$ is a homologically trivial sphere.

Theorem II is stated as Theorem 11.4.7 in [GS] and is attributed to
Lisca and Mati\'c [LM, Theorem 5.2]. 
(For spheres only the inequality $S^2\le -2$ is stated in [GS].)
The result was stated and proved in [LM] only for the case $S^2\ge 0$, 
the reason being that the paper [LM] preceeded the work of Ozsv\'ath 
and Szab\'o [OS] on surfaces with negative self-intersection. 
Joining [LM] and [OS] one immediately obtains the general case of 
Theorem II as follows. (A similar proof is given in [GS].)

By Theorem 3.2 and Corollary 3.3 in [LM] 
every relatively compact domain $\Omega\ss X$ in a Stein surface
admits a biholomorphic map onto a domain in a compact K\"ahler surface $Y$ 
with $b_2^+(Y)>1$ whose canonical bundle $K_Y =\L^2 T^*Y$ is ample, 
meaning in particular that $c_1(K_Y)^2>0$ [H, p.\ 365]. 
(Equivalently, the canonical homology class $[K_Y]= PD(c_1(K_Y)) \in H_2(X;Z)$
satisfies $[K_Y]^2>0$.) From $\L^2(TY)\simeq K_Y^*= K_Y^{-1}$ we see that
$c_1(Y)=-c_1(K_Y)$ and hence $c_1(Y)^2>0$ as well.
Such  surface is always {\it minimal} (without complex spheres 
of self-intersection number $-1$) and {\it of general type}. 
Furthermore, if $S\subset \Omega$ is homologically nontrivial in $X$ then its image 
is nontrivial in $Y$. The right hand side of (*) does not change if we replace 
$S\subset X$ with its image in $Y$. The inequality (*) now follows from Theorem I 
unless $S$ is a sphere with $S^2<0$. If (*) fails for such $S$ then by [FiS] 
(or by Theorem 1.3 in [OS]) $K=c_1(Y)+ 2\e PD(S)$ is a Seiberg-Witten basic class, 
where $\e=\pm 1$ is the sign of $c_1(Y)\cdotp S$. Since $Y$ is minimal and of general type, 
its only Seiber-Witten basic classes are $\pm c_1(Y)$ [Wi, T2], hence
$c_1(Y)+ 2\e PD(S)=\pm c_1(Y)$. This gives either $PD(S)=0$ (a contradiction) 
or $PD(S)=-\e c_1(Y)$ in which case $S^2=PD(S)^2=c_1(Y)^2>0$, a contradiction to 
the assumption $S^2<0$. 
\endpr

\ni \it Remark. \rm
Theorem II appears as Theorem 9 in the paper of S.\ Nemirovski [N]
(1999) without acknowledging the 1997 work of Lisca and Mati\'c 
[LM], even though [LM] is included in the bibliography in [N].
(In fact, all relevant references cited above were already available at that time 
and are included in [N].) Nemirovski gave in [N] a slightly different 
reduction argument to the K\"ahler case which nevertheless used essentially the same tools 
and results. He begins by embedding a domain $\Omega\subset X$ containing the given 
surface $S$ into a smooth affine algebraic surface $V\subset \C^N$ using a theorem 
of Stout [S] and Demailly, Lempert and Shiffman [DLS]. By taking the projective 
closure of $V$ in $\CP^N$ and desingularizing at infinity one obtains 
a compact K\"ahler surface $Y$ containing a biholomorphic copy of $S$. 
(So far the argument is the same as in [LM], except that $Y$
does not have the additional properties.) Nemirovski then argues that 
$S^2\ge 0\Rightarrow b_2^+(Y)>1$ and hence Theorem I applies, while $S^2<0$ 
should imply by [OS] that one of the classes $\pm [S]$ is represented by a 
complex curve, but this cannot happen by the construction of $Y$. 
However, the conclusion in Theorem 1.3 of [OS] is more complicated
when $b_2^+(Y)=1$ and it was not verified in [N] whether this theorem 
can indeed be applied as indicated. The argument given above avoids 
this situation.
\endpr

For immersed real surfaces in Stein surfaces one has the following 
adjunction inequality (see Remark 3 in [N, p.\ 742]).

\proclaim Theorem III: Let $\pi\colon S\to X$ be an immersed closed oriented real
surface with simple double points in a Stein surface $X$. If $\pi$ has $\d_+$ 
positive double points and $[\pi(S)] \ne 0$ in $X$ then 
$$  g(S)+ \d_+ 
	\ge 1+ {1\over 2}
       \left( \pi(S)^2+|c_1(X)\cdotp \pi(S)| \right).     \eqno(**) 
$$

Lacking a precise reference we include a sketch of proof which is 
intended solely for the non-experts in topology 
(no originality whatsoever is being claimed).
The first step is exactly as in the proof of Theorem II: we replace
$X$ by a compact K\"ahler minimal surface of general type, with $b_2^+(Y)>1$
and $[K_Y]^2>0$, such that $\pi(S)$ is homologically nontrivial in $Y$. From now
on we consider $\pi$ as an immersion of $S$ into $Y$. 
We then replace each positive double point of $\pi(S) \subset Y$ 
by an embedded handle as in the proof of Theorem 4.2. This increases the genus of the 
immersed surface by $\d_+$ but it does not change its homology class in $Y$, hence 
the right hand side of (**) remains unchanged. We denote the new surface again by $S$.

It remains to desingularize $\pi(S)$ by blowing up $Y$ at each of the remaining 
(negative) double points. This is described in Sec.\ 2.2 of [GS] and goes as follows.
(I thank S.\ Ivashkovich for explaining me the basic idea.)
Locally near a negative double point $p\in \pi(S)\subset Y$ we take as the local 
model for $\pi(S)$ the union $L_1\cup \bar L_2 \subset \C^2$, where $L_1 \ne L_2$ 
is a pair of complex lines through $p=0\in\C^2$ and the bar on $L_2$ means that 
the orientation has been reversed (in order to have a negative intersection at $0$). 
Let $\tau \colon \wt Y\to Y$ denote the blow-up of $Y$ at $p$ with the exceptional
sphere $e=\tau^{-1}(p)$. 
Denote by $\wt S=\tau^{-1}(\pi(S)) \subset \wt Y$ the total transform of $\pi(S)$ 
and by $S' = \bar{\tau^{-1}(\pi(S) \bs \{p\})} \subset \wt Y$ 
its proper transform. By general properties of blow-ups the right hand side 
of (**) does not change if we replace $\pi(S)\subset Y$ by the total
transform $\wt S\subset \wt Y$. We claim that the same is true if we replace 
$\wt S$ by the proper transform $S'\subset \wt Y$ of $\pi(S)$. 
This follows immediatelly from $[S']=[\wt S] \in H_2(\wt Y;Z)$ 
which can be seen as follows. We obtain $[S']$ by subtracting from 
$[\wt S]$ a copy of $[e]$ (to obtain the proper transform along $L_1$) 
as well as a copy of $-[e]$ (for the proper transform along $\bar L_2$; 
the minus sign accounts for the reversed orientation). Thus $[e]$ cancels out
and we have $[S']=[\wt S]$ as claimed. In particular, $S'\cdotp e= \wt S\cdotp e=0$.

Performing blow-up at all negative double points of $\pi(S)$ we 
get a compact K\"ahler surface $Y_0$ with $b_2^+(Y_0)>1$,
along with a natural projection $\tau \colon Y_0\to Y$,
and a smoothly embedded, oriented, homologically nontrivial real surface 
$S_0\subset Y_0$ with $\tau(S_0)=\pi(S) \subset Y$ whose genus 
$g(S_0)$ equals the genus of the original surface (immersed in $X$) 
plus $\d_+$. By Theorem I the inequality (*) holds for $S_0\subset Y_0$ 
(and hence (**) holds for the initial immersion into $X$)  
unless $S_0$ is a sphere with $S_0^2<0$ (this happens if 
we start with an immersed sphere with only negative double points). 

In the latter case we proceed as follows. Let $e_j=\tau^{-1}(p_j)$ be 
the exceptional spheres over the blown-up points and $E_j$ their 
Poincar\'e duals. Since (*) fails, we conclude as in the proof of Theorem II 
(using [FiS] or [OS]) that $c_1(Y_0)+2\e PD(S_0) \in H^2(Y_0;Z)$ is a 
Seiberg-Witten basic class for one of the choices of $\e=\pm1$. 
By [GS, Theorem 2.4.9] the only Seiberg-Witten basic classes of $Y_0$ are 
$\pm \tau^* c_1(Y)+\sum \e_j E_j$ where $\e_j=\pm 1$ for every $j$.
Also $c_1(Y_0)=\tau^* c_1(Y)-\sum E_j$. If follows that
$PD(S_0)=\pm \tau^* c_1(Y)$ modulo the exceptional classes $E_j$. 
Dualizing gives $[S_0]=\pm \tau^{-1} [K_Y]$ modulo the exceptional 
spheres $[e_j]$. (Here $[K_Y]$ denotes the Poincar\'e dual of 
$c_1(K_Y) = -c_1(Y)$.) Pushing down by $\tau\colon Y_0\to Y$ 
we get $[\pi (S)]=\pm [K_{Y}]$ and hence $[S_0]^2= [\pi(S)]^2=[K_Y]^2>0$ 
in contradiction to the assumption $S_0^2<0$.
\endpr

\ni \it Problem. \rm 
Prove Theorems II and III without using the Seiberg-Witten theory.

\demo Acknowledgements: I wish to thank S.\ Ivashkovich and M.\ Slapar for 
stimulating and helpful discussions. In particular, Ivashkovich pointed 
out to me Gromov's theorem which is used in Corollary 4.4 and explained to
me the blow-up procedure at negative double points in the proof 
of Theorem III in the Appendix. Slapar helped me substantially with 
references to gauge theory. This research was supported in 
part by a grant from the Ministry of Science and Education of the 
Republic of Slovenia.

\medskip\ni\bf References. \rm\medskip

\item{[AF]} A.\ Andreotti, T.\ Frankel:
The Lefschetz theorem on hyperplane sections. 
Ann.\ of Math.\ (2) {\bf 69} (1959), 713--717.

\item{[B]} E.\ Bishop: 
Differentiable manifolds in complex Euclidean spaces.
Duke Math.\ J.\ {\bf 32} (1965), 1--21.

\item{[BK]} E.\ Bedford and W.\ Klingenberg:
On the envelope of holomorphy of a two-sphere in $\C^2$.
J.\ Amer.\ Math.\ Soc.\ 4 (1991), 623--646.

\item{[CS]} S.\ Chern and E.\ Spanier:
A theorem on orientable surfaces in four-dimen\-sional space.
Comm.\ Math.\ Helv.\ {\bf 25} (1951), 205--209.

\item{[DLS]} J.-P.\ Demailly, L.\ Lempert and B.\ Schiffman:
Algebraic approximations of holomorphic maps from Stein domains
to projective manifolds.
Duke Math.\ J.\ {\bf 76} (1994), 333--363.

\item{[E1]} Y.\ Eliashberg: 
Topological characterization of Stein manifolds of dimension $>2$. 
Internat.\ J.\ Math. {\bf 1} (1990), 29--46. 

\item{[E2]} Y.\ Eliashberg: Filling by holomorphic discs and 
its applications. 
Geometry of low-dimensional manifolds, 2 (Durham, 1989), 45--67, 
London Math.\ Soc.\ Lecture Note Ser., 151, 
Cambridge Univ.\ Press, Cambridge, 1990. 

\item{[FiS]} R.\ Fintushel and R.\ Stern: 
Immersed spheres in $4$-manifolds and the immersed Thom conjecture.
Turk.\ J.\ Math.\ {\bf 19} (1995), 145--157.

\item{[F1]}  F.\ Forstneri\v c: 
Analytic discs with boundaries in a maximal real submanifold of $\C^2$.
Ann.\ Inst.\ Fourier {\bf 37} (1987), 1--44.

\item{[F2]}  F.\ Forstneri\v c: 
Complex tangents of real surfaces in complex surfaces. 
Duke Math.\ J.\ {\bf 67} (1992), 353--376. 

\item{[FSt]} F.\ Forstneri\v c, E.\ L.\ Stout: 
A new class of polynomially convex sets. 
Ark.\ Mat.\ {\bf 29} (1991), 51--62.

\item{[Go]} R.\ E.\ Gompf: Handlebody construction of Stein surfaces. 
Ann.\ of Math.\ (2) {\bf 148} (1998), 619--693. 

\item{[GS]} R.\ E.\ Gompf, A.\ I.\ Stipsicz: 
{\it $4$-manifolds and Kirby Calculus}. 
American Mathematical Society, Providence,  1999.

\item{[Gra]} H.\ Grauert: On Levi's problem and the 
imbedding of real-analytic manifolds. 
Ann.\ of Math.\ (2) {\bf 68} (1958), 460--472.

\item{[Gr]} M.\ Gromov: Partial Differential Relations.
Ergebnisse der Mathematik und ihrer Grenzgebiete (3) {\bf 9},
Springer, Berlin--New York, 1986.

\item{[GH]} P.\ Griffiths, J.\ Harris: 
{\it Principles of Algebraic Geometry.} 
Pure and Applied Mathematics, Wiley-Interscience, New York, 1978.

\item{[H]} R.\ Hartshorne: {\it Algebraic Geometry}.
Graduate Texts in Mathematics {\bf 52}, Springer, Berlin, 1977.

\item{[IS]} S.\ Ivashkovich and V.\ Shevchishin:
Structure of the moduli space in a \nbd\ of a cusp-curve 
and meromorphic hulls.
Invent.\ Math.\ {\bf 136} (1999), 571--602.

\item{[J]} B.\ J\"oricke: Local polynomial hulls of discs near 
isolated parabolic points.
Indiana Univ.\ Math.\ J.\ {\bf 46} (1997), 789--826.

\item{[KW]} C.\ E.\ K\"onig, S.\ M.\ Webster:
The local hull of holomorphy of a surface in the space of two complex variables.
Invent.\ Math.\ {\bf 67} (1982), 1--21.

\item{[KM]} P.\ B.\ Kronheimer, T.\ S.\ Mrowka:
The genus of embedded surfaces in the projective plane. 
Math.\ Res.\ Lett.\ {\bf 1} (1994), 797--808. 

\item{[Kr]}  N.\ G.\ Kruzhilin: Two-dimensional spheres in the
boundaries of strictly pseudoconvex domains in $\C^2$.
Izv.\ Akad.\ Nauk SSSR Ser.\ Mat.\ {\bf 55} (1991), 1194--1237.
(English transl.\ in Math.\ USSR Izv.\ {\bf 39} (1992), 1151--1187.) 

\item{[La]} H.\ F.\ Lai: Characteristic classes of real manifolds 
immersed in complex manifolds. 
Trans.\ Amer.\ Math.\ Soc.\ {\bf 172} (1972) 1--33.

\item{[LM]} P.\ Lisca, G.\ Mati\'c: 
Tight contact structures and Seiberg-Witten invariants. 
Invent.\ Math.\ {\bf 129} (1997), 509--525. 

\item{[M]} J.\ Morgan: {\it The Seiberg-Witten Equations
and Applications to the Topology of Smooth Four-Manifolds.}
Math.\ Notes {\bf 44}, Princeton University Press,
Princeton, 1996.

\item{[MST]}  J.\ W.\ Morgan, Z.\ Szab\'o, C.\ H.\ Taubes:
A product formula for the Seiberg-Witten invariants and 
the generalized Thom conjecture. 
J.\ Differential\ Geom.\ {\bf 44} (1996), 706--788. 

\item{[N]} S.\ Nemirovski: 
Complex analysis and differential topology on complex surfaces.
Uspekhi Mat.\ Nauk {\bf 54}, 4 (1999), 47--74.
(English transl.\ in Russian Math.\ Surveys {\bf 54}, 4 (1999), 729--752.)

\item{[OS]} P.\ Ozsv\'ath and Z.\ Szab\'o: 
The symplectic Thom conjecture.
Ann.\ of Math.\ (2) {\bf 151} (2000), 93--124. 

\item{[R]} W.\ Rudin: A totally real Klein bottle in $\C^2$. 
Proc.\ Amer.\ Math.\ Soc.\ {\bf 82} (1981), 653--654.

\item{[S]} E.\ L.\ Stout: Algebraic domains in Stein manifolds.
Proceedings of the conference on Banach algebras and several complex 
variables (New Haven, Conn., 1983), 259--266, 
Contemp.\ Math.\ {\bf 32}, Amer.\ Math.\ Soc., Providence, RI, 1984.

\item{[T1]} C.\ Taubes: The Seiberg-Witten invariants and symplectic forms.
Math.\ Res.\ Lett.\ {\bf 1} (1994), 809--822.

\item{[T2]} C.\ Taubes: 
SW$\Rightarrow$Gr: from the Seiberg-Witten equations to pseudo-holo\-mor\-phic curves.
J.\ Amer.\ Math.\ Soc.\ {\bf 9} (1996), 845--918.

\item{[We]} S.\ M.\ Webster:
Minimal surfaces in a K\"ahler surface.
J.\ Diff.\ Geom.\ {\bf 20} (1984), 463--470. 

\item{[Wn]} A\ Weinstein: {\it Lectures on Symplectic Manifolds.}
Reg.\ Conf.\ Ser.\ Math.\ {\bf 29}, Amer.\ Math.\ Soc., Providence, 1977.

\item{[Wg]} J.\ Wiegerinck: 
Local polynomially convex hulls at degenerated CR singularities 
of surfaces in $\C^2$. 
Indiana Univ.\ Math.\ J.\ {\bf 44} (1995), 897--915.

\item{[Wi]} E.\ Witten:  Monopoles and four-manifolds.
Math.\ Res.\ Lett.\ {\bf 1} (1994), 769--796.

\bigskip\medskip
\ni Institute of Mathematics, Physics and Mechanics, 
University of Ljub\-ljana,
Jadranska 19, 1000 Ljubljana, Slovenia

\bye